\documentclass[notitlepage,leqno,12pt]{article}
\usepackage{amssymb}
\catcode`\@=11 \@addtoreset{equation}{section}

\catcode`\@=12
\usepackage[utf8]{inputenc}
\usepackage{latexsym}
\usepackage{textcomp}
\usepackage{amsmath}
\usepackage{amsfonts}
\usepackage{amssymb}
\usepackage{mathrsfs}
\usepackage{times}
\renewcommand{\a }{\alpha}

\renewcommand{\d}{\delta}

\newcommand{\e }{\varepsilon }

\renewcommand{\l }{\lambda }

\newcommand{\s }{\sigma }

\renewcommand{\o }{\omega }
\renewcommand{\O }{\Omega }

\newcommand{\be}{\begin{equation}}
\newcommand{\ee}{\end{equation}}

\newcommand{\R}{\mathbb{R}}
\newcommand{\N}{\mathbb{N}}

\newcommand{\C}{\mathbb{C}}

\newcommand{\M}{\mathcal{M}}

\newtheorem{Theorem}{Theorem}[section]
\newtheorem{Lemma}[Theorem]{Lemma}

\newtheorem{Corollary}[Theorem]{Corollary}

\def\proof{\noindent{{\bf Proof. }}}
\def\square{\vbox{
\hrule height .4pt
\hbox{\vrule width .4pt height 7pt \kern 7pt
\vrule width .4pt}
\hrule height .4pt }}
\def\QED{\hfill {$\square$}\goodbreak \medskip}
\def\R{{\mathbb R}}

\def\C{{\mathcal C}}

\font\sc=cmcsc9  \textwidth=14truecm
\author{El hadji Abdoulaye Thiam
\footnote{\footnotesize{AIMS-Senegal.
E-mail: \textit{elhadji@aims-senegal.org},
\textit{heat1719@gmail.com}. } }}
\begin{document}
\title{WEIGHTED HARDY INEQUALITY ON RIEMANNIAN MANIFOLDS}
\date{}
\maketitle
\bigskip
\noindent{\footnotesize{\bf Abstract.}
Let $(\M,g)$ be a smooth compact Riemannian manifold of dimension $N\geq 3$ and we let $\Sigma$ to be a closed submanifold of dimension $1 \leq k \leq N-2. $  In this paper we study existence and non-existence of minimizers of Hardy inequality with weight function singular on $\Sigma$ within the framework of Brezis-Marcus-Shafrir \cite{BMS}. In particular we provide necessary and sufficient conditions for existence of minimizers.
\section{Introduction}\label{section 1}
Let $N \geq 3$ , $1 \leq k \leq N-2$ and pose $x=(y,z)\in \R^k \times \R^{N-k}$. We denote by $\mathcal{D}^{1,2}(\R^N)$ the completion of $\mathcal{C}^\infty_c(\R^N)$ with respect to the norm
$$
u\longmapsto \int_{\R^N} |\nabla u|^2 dx.
$$  
We recall the following Hardy type inequality with cylindrical weights
\begin{equation}\label{A2}
\int_{\R^N} |\nabla u|^2 dx \geq \biggl(\frac{N-k-2}{2}\biggl)^2 \int_{\R^N}|y|^{-2} |u|^2 dx, \;\ \forall u \in \mathcal{D}^{1,2}(\R^N),
\end{equation}
see the book of May'ja \cite{Mayza} for a proof. See also the work of Brezis-Vasquez \cite{BV}, Musina \cite{Musina} and Gazzini-Musina \cite{GM}.
The constant $\bigl(\frac{N-k-2}{2}\bigl)^2$ is sharp and never achieved in $\mathcal{D}^{1,2}(\R^N)$. However there exists a function
$$u(x)=|y|^{\frac{2+k-N}{2}}$$
which satisfies
\begin{equation}
\Delta_{\R^N} u +\biggl(\frac{N-k-2}{2}\biggl)^2 |y|^{-2} u =0.
\end{equation}
We will call this function a "virtual" ground state because it does not belong to $\mathcal{D}^{1,2}(\R^N)$.
Inequality \eqref{A2} is invariant by translation in the $z$-variable and by scaling in the full variable  yielding difficulties in the study of elliptic and parabolic equations involving inverse square potentials.\\\
Note that \eqref{A2} is not in general valid for a Riemannian manifold ($\mathcal{M}^N,g$) of dimension $N\geq 3$. However, by Allegretto-Piepenbrink argument (see \cite{Allegretto} and \cite{Piepenbrink}) and by construction of super-solution near $\Sigma,$ we prove the local Hardy inequality below in Lemma \ref{lemma31}, in a small tubular neighborhood 
$$\Sigma_r:=\lbrace p \in \mathcal{M}: \rho(p):= dist(p, \Sigma)< r\rbrace$$
of $\Sigma$, i.e
\begin{equation}\label{Local}
\int_{\Sigma_r} |\nabla u| dv_g \geq \biggl(\frac{N-k-2}{2}\biggl)^2 \int_{\Sigma_r} \rho^{-2} |u|^2 dv_g, \;\ \forall u \in H^1({\Sigma_r}),
\end{equation}
where $\rho(p):=dist(p, \Sigma)$ is the geodesic distance to $\Sigma.$
This type of result was first proved by Brezis-Marcus in \cite{BM}. See also the work of Fall-Mahmoudi in \cite{FMa} and  Thiam in \cite{THIAM}.\\\
Using \eqref{Local} with an argument of partition of unity around $\Sigma$, we will prove in Lemma \ref{lemma4} below the following
\begin{equation}\label{1.4}
\int_{\M} |\nabla u|_g^2 dv_g \geq \biggl(\frac{N-k-2}{2}\biggl)^2 \int_\M \rho^{-2} |u|^{2} dv_g + \l \int_\M u^2 dv_g, \;\ \forall u \in H^1(\M)
\end{equation}
for a constant $\l$ depending on $\M$.
We remark that the Hardy inequality is a particular case of the Caffarelli-Kohn-Nirenberg inequality, see \cite{CKN}. The knowledge of Hardy, Hardy-Sobolev, Gagliardo-Nirenberg, Sobolev or Caffarelli-Kohn-Nirenberg inequality on a manifold $\M$ and their best constants allows to obtain qualitative properties on the manifold $\M$. For instance in \cite{AX}, \cite{CX} and \cite{XIA} it was shown that if $\M$ is a complete open Riemannian manifold with non negative Ricci curvature in which a Hardy or Gagliardo-Nirenberg or Caffarelli-Kohn-Nirenberg type inequality holds, then $\M$ is in some suitable sense close to the Euclidean space.\\\
Inequalities involving integrals of a function and its derivatives together with singular weights appear frequently in various branches of mathematics and represent a useful tool in the theory of differential equations. They have several applications in many questions from mathematical physics, spectral theory, analysis of linear and nonlinear PDEs, harmonic analysis and stochastic analysis.
For more details related to these inequalities, in particular the Hardy one, see \cite{AX, BMR, BFT, BM, Caron, Amb1, Amb2, DA, FMS, Fall, FMa, FTT, GGM, LL1, KMP, KO, ME, Markus, Mats}.\\\
In this paper, we are interested in the following Hardy inequality with weight functions on a compact Riemannian manifold $(\M,g)$ of dimension N. Therefore we propose to study the problem of finding minimizers of the following quotient in the spirit of Brezis-Marcus \cite{BM}
\begin{equation}\label{Principal}
\mu_\l(\M,\Sigma,b,q,\eta):=\inf_{u \in H^1(M)} \displaystyle\frac{\displaystyle\int_\M b |\nabla u|_g^2 dv_g- \lambda \int_\M \rho^{-2} |u|^2 \eta dv_g}{\displaystyle\int_\M \rho^{-2} |u|^2 q dv_g} ,
\end{equation}
where $\rho(p):=\textrm{dist}(p,\Sigma)$ is the geodesic distance function to $\Sigma$ and the weights functions $b,\,\ q $ and $\eta$ satisfy
\begin{equation}\label{A22}
b, q \in \mathcal{C}^2(\M),\;\;\ b, q >0 \;\;\ \textrm{in}\;\;\ \M,\;\;\ \eta>0 \;\;\ \textrm{in} \;\;\ \M\setminus \displaystyle\Sigma,\;\;\ \eta \in Lip(\M)
\end{equation}
and
\begin{equation}\label{A33}
\max_{\displaystyle\Sigma} \frac{q}{b}=1,\;\;\;\;\;\;\;\;\;\;\ \eta=0 \;\;\;\ \textrm{on} \;\;\ \displaystyle\Sigma.
\end{equation}
We have the following
\begin{Theorem}\label{Theorem 1}
Let $(\M,g)$ be a smooth compact Riemannian manifold of dimension $N \geq 3$ and let $\Sigma\subset \M$ be a closed submanifold of dimension $k=1,...,N-2$. Assume that the weight functions $b,q$ and $\eta$ satisfy \eqref{A22} and \eqref{A33}.
Then, there exists $\lambda^*=\lambda^*(b,q,\eta,\M,\Sigma)$ such that
\begin{equation}\label{EL1}
\displaystyle \mu_\l(\M,\Sigma)= \biggl( \displaystyle\frac{N-k-2}{2}\biggl)^2, \;\;\ \forall \lambda \leq \lambda^*,
\end{equation}
\begin{equation}\label{EL2}
\displaystyle \mu_\l(\M,\Sigma)< \biggl( \displaystyle\frac{N-k-2}{2}\biggl)^2, \;\;\ \forall \lambda > \lambda^*.
\end{equation}
The infinimum $\mu_\l(\M,\Sigma)$ is attained if $\l>\l^*$ and it is not attained when $\l<\l^*$.
\end{Theorem}
The existence of $\l^*$ is a consequence of the local Hardy inequality 
\begin{equation}
\int_{\Sigma_r} b |\nabla u|^2 dv_g \geq \biggl(\frac{N-k-2}{2}\biggl)^2 \int_{\Sigma_r} q |u|^2 \rho^{-2} dv_g
\end{equation}
(see Lemma \ref{lemma31} and Lemma \ref{lemma4}). The existence and non-existence parts are classic. They were almost the same done in \cite{BM} and in \cite{THIAM}. A natural question is to know what happens concerning the critical case. Thus, we have the following
\begin{Theorem}\label{Theorem 2}
let $\l^*$ be given by Theorem \ref{Theorem 1}.
Then $\mu_{\l^*}(\M, \Sigma)$ is achieved if and only if 
\begin{equation}\label{Terence}
\int_{\Sigma} \frac{d\s}{\sqrt{1-q(\s)/b(\s)}}< \infty.
\end{equation}
\end{Theorem}
As a consequence of this, we get the following
\begin{Corollary}\label{Corollary}
Let $(\mathcal{M},g)$ be a compact Riemannian manifold of dimension $N \geq 3$ and let $\Sigma$ be a closed submanifold of dimension $1\leq k \leq N-2$. For $\l \in \R$, put
\begin{equation}
\mathcal{V}_\l(\mathcal{M},\Sigma):= \inf_{u \in H^1(\mathcal{M})} \displaystyle \frac{\displaystyle\int_{\mathcal{M}} |\nabla u|^2 dv_g-\l \int_{\mathcal{M}} u^2 dv_g}{\displaystyle\int_{\mathcal{M}} \rho^{-2} u^2 dv_g}.
\end{equation}
Then there exists $\l^*=\l^*(\mathcal{M}, \Sigma)$ such that
$$\mathcal{V}_\l(\mathcal{M},\Sigma)=\displaystyle \biggl(\frac{N-k-2}{2}\biggl)^2, \;\ \forall \l \leq \l^*,$$
$$\mathcal{V}_\l(\mathcal{M},\Sigma)< \displaystyle \biggl(\frac{N-k-2}{2}\biggl)^2, \;\ \forall \l > \l^* .$$
Moreover $\mathcal{V}_\l(\mathcal{M},\Sigma)$ is attained if and only if $\l>\l^*.$ 
\end{Corollary}
When the singularity is reduced to a single point $\lbrace p_0\rbrace$ ($k=0$), the corollary remain valid. It was proved by Thiam in \cite{THIAM}.\\
Our arguments of proof are based on the construction of a $H^1$ super-solution and a $H^1$ sub-solution of the linear operator $L_\l$ defined in \eqref{lll}. Without any loss of generality we may assume that $b\equiv 1$ (see Section \ref{Section 4} below) and $r$ is small enough and we perturb the virtual ground-state $$v_{a,q}(p)=(-log \rho)^a \rho^\a(p)$$ for the Hardy constant $\bigl(\frac{N-k-2}{2}\bigl)^2$, where 
$$\displaystyle \alpha(x)=\frac{2+k-N}{2}\biggl(1- \sqrt{1- q(\s(x))+|x|}\biggl).$$
Furthermore, it's easy to verify that for $a < - \frac{1}{2}$ and for $\e \in (0,1)$, $v_{a,q}$ and $v_{0,q-\e}$ belong to $H^1(\Sigma_r)$.
We prove the non-existence part by assuming by contradiction that, when 
\begin{equation}\label{123456}
\int_{\Sigma} \frac{d\s}{\sqrt{1-q(\s)}}=\infty,
\end{equation}
there exists a non-negative solution $u \in H^1(\M) \cap \C(\M\setminus \Sigma).$ We then construct a $H^1(\Sigma_r)$ sub-solution $V_\e:=v_{-1,q}+v_{0, q-\e}$ which is upper bounded by $u$ ( modulo a multiplicative positive constant independent on $a$ and $\e$) so that 
\begin{equation}\label{B71}
||\rho^{-1}V_\e||_{L^2(\Sigma_r)} \leq C ||\rho^{-1}u||_{L^2(\M)} \leq C' ||u||_{H^1(\M)}
\end{equation}
by the Hardy inequality \eqref{1.4}.
Moreover using polar coordinates we  verify that
\begin{equation}\label{B72}
C \int_{\Sigma} \frac{d \s}{\sqrt{1-q(\s)}}\leq \int_{\Sigma_r} V_0^2 \rho^{-2} dv_g
\end{equation}
for $r$ small enough.
Hence taking the limit in \eqref{B71} as $\e \longrightarrow 0$, we get contradiction.\\\
For the existence part, we construct a super-solution $U:= v_{0,q}-v_{-1,q}$ and we suppose that
\begin{equation}\label{HO2}
\int_{\Sigma} \frac{d \s}{\sqrt{1-q(\s)}} < +\infty.
\end{equation}
Then $U \in H^1(\Sigma_r),$ (see Lemma \ref{lemma2} below). Next, we let the sequence of real numbers $\lbrace \l_n\rbrace$ decreasing to $\l^*$. By Theorem \ref{Theorem 1}, we can now associate to each $\l_n$ a positive minimizer $u_n \in H^1(\M) \cap \C(\M\setminus \Sigma)$ for $\mu_{\l_n}$. Then using some comparison argument, the sequence $\lbrace u_n \rbrace$ is uniformly bounded in $\Sigma_{r_0}$ by the super-solution $U$ (modulo a multiplicative positive constant independent on $n$). Hence $\rho^{-1} u_n$ converge strongly to $\rho^{-1}u$ in $L^2(\M)$ by Rellich-Kondrakov theorem and that $ u_n$ converge to $u$ in $H^1(\M)$ strongly.\\
The paper is organized as follows. In Section \ref{section 2}, we give some Preliminaries and Notations and we construct a super and a sub-solutions we will use in Section \ref{Section 4} to prove Theorem \ref{Theorem 2} and in Section \ref{section 3}, we prove the existence of $\l^*$ and we give a complete proof of Theorem \ref{Theorem 1}.\\
\textbf{ACKNOWLEDGEMENTS}\\
I wish to thank my supervisor Mouhamed Moustapha Fall for his disponibility. This work is supported by the German Academic Exchange Service(DAAD).
\section{Preliminaries and Notations}\label{section 2}
Consider $p \in \Sigma$. We denote by $T_p \Sigma$ the tangent space of $\Sigma$ and $N_p \Sigma$ he normal space of $T_p \Sigma$ at $p$.  We may assume that
\begin{equation}
N_p\Sigma=\textrm{Span} \bigl< E_1,..., E_{N-k}\bigl> \;\;\;\ \textrm{and} \;\;\;\ T_p\Sigma= \textrm{Span} \bigl<E_{N-k+1},..., E_N\bigl>.
\end{equation}
A neighborhood of $p$ in $\Sigma$ can be parametrized via the mapping
\begin{equation}\label{Sur1}
\textrm{Exp}_p^{\Sigma}: B_r(0) \subset \R^{k}\rightarrow \Sigma_r \supset \Sigma
\end{equation}
$$
\begin{array}{ll}
y \longmapsto f^p(z)= \textrm{Exp}_p^{\Sigma}\biggl(\sum_{a=N-k+1}^N z_a E_a \biggl),
\end{array}
$$
where $z=(z_{N-k+1},...,z_{N}) \in \R^{k}$, $B_r(0)$ is the ball centered at $0$ and of radius $r$,  $ \textrm{Exp}_p^{\Sigma} $ is the exponential mapping at $p$ in $\Sigma$ and $\Sigma_r$ defined in \eqref{KKK}.
Now we extend $(E_i)_{1 \leq i\leq N-k}$ to an orthonormal frame $(X_i)_{1\leq i\leq N-k}$ in a neighborhood of $p$ in $\M$ via the mapping
\begin{equation}\label{Sur2}
\textrm{Exp}_{f_p(z)}^{\M}:\R^k \times \R^{N-k}\rightarrow \M
\end{equation}
$$
\begin{array}{ll}
x=(y,z) \longmapsto F_{\M}^{p_i}(x)=\textrm{Exp}_{f_p(z)}^{\M} \biggl(\sum_{i=1}^{N-k} y_i X_i \biggl),
\end{array}
$$
where $y=(y_1,...,y_{N-k})$ and $\textrm{Exp}_{f_p(z)}^{\M}$ is the exponential map at $f_p(z)$ in $\M$.\\
In the following, we will consider the geodesic neighborhood contained in $\M$ around $\Sigma$ of radius $r$
\begin{equation}\label{KKK}
\Sigma_r= \bigl \lbrace p \in \M: \rho(p):=dist(p, \Sigma) < r\bigl\rbrace.
\end{equation}
In these normal coordinates, the Laplace-Beltrami operator is given by
\begin{equation}\label{LBO}
\Delta_g=- g^{ij}\biggl(\frac{\partial^2}{\partial_{x_i} \partial_{x_j}} - \Gamma^k_{ij} \frac{\partial}{\partial_{x_k}}\biggl),
\end{equation}
where $\bigl\lbrace \Gamma^k_{ij} \bigl\rbrace_{1\leq i,j,k\leq N}$ are the components of the metric $g$ and $g^{ij}=(-g^{-1})_{ij}$ are the components of the inverse matrix of $g$. Then the following estimates hold
\begin{equation}\label{Estimation2}
\Gamma_{ij}^k(x)=O(|y|), \;\;\;\;\;\   g_{ij}(F_\M^{p_i}(x))= \delta_{ij}+ O(|y|^2) \;\;\;\ \textrm{and} \;\;\;\;\;\ \rho_k(F_\M^{p_i}(x))= |y|,
\end{equation}
see the paper of Mahmoudi-Mazzeo-Pacard \cite{MMP}.
In addition, there exists a positive constant $r_0$ depending on $\Sigma$ and $\M$ such that $\rho \in \mathcal{C}^{\infty}_c(\Sigma_r)$.
Moreover $\Sigma$ is a closed submanifold of a compact manifold $\M$, then for $r$ sufficiently small, there exists a finite number of Lipschitz open sets $(\Omega_i)_{1\leq i \leq  N_0}$ such that
$$
\O_i \cap \O_j=\emptyset \;\;\ \textrm{for} \;\;\;\ i \neq j
$$
and
\begin{equation}\label{A5}
\overline{\Sigma}_r= \displaystyle\bigsqcup_{i=1}^{N_0} \overline{\O}_i.
\end{equation}
We choose the open sets $\O_i$, using the above Fermi coordinates, so that
\begin{equation}
\O_i=F_M^{p_i}\bigl(B^{N-k}(0,r) \times D_i\bigl) \;\;\;\ \textrm{with} \;\;\ p_i \in \Sigma,
\end{equation}
where the $D_i$'s are Lipschitz disjoint open sets of $\R^k$ such that 
\begin{equation}
\bigcup_{i=1}^{N_0} \overline{f^{p_i}(D_i)}= \Sigma.
\end{equation} For $p \in \M$, we denote by $\s(p)$ the orthogonal projection of $p$ on $\Sigma.$ For the rest of the paper, if there is no confusion, we use the notation $v_a$ instead of $v_{a,q}$.
We get the following
\begin{Lemma}\label{lem3}
Let $a \in \R$ and define
\begin{equation}\label{label11}
v_{a,q}(p)=\bigl(-\textrm{log} \rho(p)\bigl)^a \rho(p)^{\alpha}
\end{equation}
where for $x=F^{-1}(p) \in \R^N$
\begin{equation}\label{alpha}
\displaystyle \alpha(x)=\frac{2+k-N}{2}\biggl(1- \sqrt{1- q(\s(x))+|x|}\biggl).
\end{equation}
Then we have
\begin{equation}
\Delta_g v_{a,q}=-\bigl(\frac{N-k-2}{2}\bigl)^2 q \rho^{-2} v_{a,q}+ a(a-1) \rho^{-2}\bigl(\textrm{log} \rho\bigl)^{-2} v_{a,q}
\end{equation}
$$+\bigl(N-k-a\bigl) \rho^{-2}\bigl(\textrm{log} \rho\bigl)^{-1} v_{a,q}+ O\bigl(\textrm{log} \rho \rho^{-3/2}(-\textrm{log}\rho)^a \rho^\a\bigl) \qquad \textrm{ in }\Sigma_r.$$
\end{Lemma}
\proof
If there is no ambiguity, we will write $\o_a$ and $v_a$ instead of $\o_{a,q}$ and $v_{a,q}$, where
\begin{equation}\label{Moustapha}
X_a(x)=\bigl(-log|x|\bigl)^a, \qquad \omega(x)=|x|^{\alpha(x)} \quad \textrm{and } \quad\omega_a=X_a \omega.
\end{equation}
We can verify easily that
\begin{equation}\label{Equation1}
\Delta_{\R^N} w_a= X_a \Delta_{\R^N} \omega+ 2 \nabla X_a \nabla \omega + \omega \Delta_{\R^N} X_a.
\end{equation}
We are going to calculate term by term the expression \eqref{Equation1} using simple calculations.\\
We have that 
$$\Delta \omega =\Delta (\varphi \circ u(x)),$$
where $\varphi(t)=e^t$ and $$u(x)=\alpha(x) log(|x|)=log w.$$
But $$\Delta (\varphi\circ u(x))=\varphi''(u(x)) |\nabla u(x)|^2 + \varphi'(u(x)) \Delta u(x)$$ and $$\varphi(u(x))=\varphi'(u(x))$$ so that
\begin{equation}\label{13}
\Delta \omega=\omega\biggl[|\nabla \textrm{log} \omega|^2+ \Delta \textrm{log}\omega\biggl].
\end{equation}
Since $\textrm{log} w=\alpha(x) (\textrm{log}|x|)$, we have that 
\begin{equation}\label{14}
\Delta \textrm{log} w= \alpha \Delta log |x| + 2 \nabla \alpha \nabla (\textrm{log} |x|) + \textrm{log} |x| \Delta \alpha.
\end{equation}
Using \eqref{13}, we get 
\begin{equation}\label{15}
\Delta \alpha(x)=\alpha\biggl[\frac{1}{2} \Delta \textrm{log}\bigl(1-q(\s(x))+|x|\bigl)+ \frac{1}{4} |\nabla \textrm{log}\bigl(1-q(\s(x))+ |x|\bigl)|^2\biggl].
\end{equation}
But
$$
\begin{array}{ll}
\nabla \biggl(\textrm{log} (1-q(\s(x))+|x|\biggl)=\displaystyle\frac{-\nabla q(\s(x))+ \nabla |x|}{1-q(\s(x))+|x|}
\end{array} 
$$
and
$$
\begin{array}{ll}
\Delta \textrm{log}(1-q(\s(x))+|x|)&=\displaystyle\frac{\Delta (1-q(\s(x))+|x|)}{1-q(\s(x))+|x|}-\displaystyle\frac{|\nabla \biggl(\textrm{log} (1-q(\s(x))+|x|\biggl)|^2}{(1-q(\s(x))+|x|)^2}\\\\\
&=\displaystyle \frac{\Delta(q\circ \s(x))+\Delta |x|}{1-q\circ \s(x)+ |x|}-\displaystyle \frac{|\nabla \biggl(\textrm{log} (1-q(\s(x))+|x|\biggl)|^2}{(1-q(\s(x))+|x|)^2}\\\\\
&=\displaystyle \frac{\Delta(q\circ \s(x))+\Delta |x|}{1-q\circ \s(x)+ |x|}- \displaystyle\frac{|\nabla q \circ \s(x)|^2 +1-2\nabla |x| \nabla (q\circ \s(x))}{\bigl( 1-q(\s(x))+|x|\bigl)^2}.\\\
\end{array}
$$
Puting the above in \eqref{15}, we obtain that
\begin{equation}\label{16}
\Delta \alpha= \alpha \biggl[\frac{1}{2} \;\;\ \frac{-\Delta (q \circ \s(x))+ \Delta |x|}{1-q\circ\s(x)+|x|}-\frac{1}{2} \;\ \frac{|\nabla q\circ \s(x)|^2+1-2\nabla |x| \nabla (q\circ\s(x))}{\bigl(1-q\circ \s(x)+ |x|\bigl)^2}
\end{equation}
$$ +\frac{1}{4} \;\ \frac{|\nabla (q\circ \s(x))|^2 +1 -2\nabla |x| \nabla q\circ \s (x)}{\bigl(1-q\circ \s(x)+ |x|\bigl)^2}\biggl].$$
Using the fact that $q \in \mathcal{C}^2,$ we conclude that
\begin{equation}\label{17}
\Delta \alpha(x)= O\bigl(|x|^{-3/2}\bigl).
\end{equation}
We have also that
$$
\begin{array}{ll}
\nabla \alpha &=\biggl(\displaystyle\frac{N-k-2}{2} \nabla \sqrt{1-q\circ \s(x) +|x|}\biggl)\\\\\
&=\displaystyle\frac{N-k-2}{2}\displaystyle\frac{1}{2 \sqrt{1-q\circ \s(x) +|x|}} \nabla (1-q\circ \s(x) +|x|).
\end{array}
$$
Therefore
$$\nabla \alpha \nabla |x|= \frac{N-k-2}{4 \sqrt{1-q\circ \s(x) +|x|}}\bigl(1-\nabla |x| \nabla (q\circ \s(x))\bigl).$$
Hence $$ \nabla \alpha \nabla |x|= O(|x|^{-1/2})$$ and from which we deduce that
\begin{equation}\label{18}
\nabla \alpha \nabla (\textrm{log} |x|)=\nabla \alpha \frac{\nabla |x|}{|x|}=O(|x|^{-3/2}).
\end{equation}
Now let us evaluate the term $ \Delta (\textrm{log}\omega)$.
We have that 
\begin{equation}
\Delta (\textrm{log} |x|)= \displaystyle\frac{N-k-2}{|x|^2}
\end{equation}
 so that 
$$ \alpha \Delta \textrm{log}|x|=\alpha \frac{N-k-2}{|x|^2}.$$
Recall that
$$\Delta( \textrm{log} \omega)= \textrm{log}(|x|)\Delta \alpha +2 \nabla \alpha \nabla (\textrm{log}|x|) +\alpha \Delta (\textrm{log}(|x|)),$$
therefore
$$\Delta \textrm{log} w=\alpha \frac{N-k-2}{2} +O (|x|^{-3/2})+O (|x|^{-3/2} \textrm{log}|x|)$$
and that 
 $$\Delta \textrm{log} w= \alpha \frac{N-k-2}{|x|^2}\bigl(1+O(|x|)\bigl) + O\bigl(|x|^{-3/2}\bigl).$$
 We have also that 
 $$\nabla (\textrm{log} w)= \nabla (\alpha \textrm{log}|x|)= \alpha \frac{\nabla |x|}{|x|}+ \textrm{log} |x| \nabla \alpha$$ and thus 
 $$ |\nabla \textrm{log} w|^2= \frac{\alpha^2}{|x|^2}+ (log |x|)^2 |\nabla \alpha|^2 +\frac{2 \alpha }{|x|} \textrm{log} |x| \nabla |x| \nabla \alpha= \frac{\alpha^2}{|x|^2}+ O\bigl(\textrm{log}|x| |x|^{-3/2}\bigl).$$
 Therefore \eqref{13} becomes
 \begin{equation}\label{aa}
 \frac{\Delta \o}{\o}= \alpha \frac{N-k-2}{|x|^2} + \frac{\alpha^2}{|x|^2} + O\bigl(\textrm{log}|x| |x|^{-3/2}\bigl).
 \end{equation}
We recall that, we want to calculate the Laplacian of $\o_a$ define in \eqref{Moustapha}. Then we have
 $$\Delta \omega_a= \o \Delta X_a  +2 \nabla X_a \nabla \omega+X_a \Delta \omega.$$ 
 From \eqref{aa} we have that
 \begin{equation}\label{nice1}
 X_a \Delta \o=\o_a \biggl[\frac{N-k-2}{|x|^2}+\frac{\alpha^2}{|x|^2}+ O\bigl(\textrm{log}|x| |x|^{-3/2}\bigl)\biggl].
 \end{equation}
 Now we are going to evaluate $\o \Delta X_a$.
 We have $$\Delta X_a=\Delta (\varphi \circ u(x)),$$
 where $\varphi(t)=(-log t)^a$ and $u(x)=|x|$. 
It's easy to verify that
 $$\Delta X_a =\varphi''(u(x)) + \displaystyle\frac{N-k-1}{|x|}.$$
Therefore
\begin{equation}\label{nice2}
\o \Delta X_a= \o_a \biggl[ \displaystyle\frac{a(a-1)}{|x|^2 (log |x|)^2} + \displaystyle\frac{N-k-2}{|x|^2 (\textrm{log} |x|)}a\biggl].
\end{equation} 
Now let us finish this part by calculate the expression $2 \nabla X_a \nabla \o.$
By simple calculations we get that
\begin{equation}\label{+}
\begin{array}{ll}
\nabla X_a = \biggl(\displaystyle\frac{a \nabla |x|}{|x| \textrm{log} |x|}\biggl) X_a
\end{array}
 \end{equation}
and
\begin{equation}\label{++}
\nabla \o= \o \biggl[\textrm{log} |x| \nabla \alpha + \alpha \frac{\nabla |x|}{|x|}\biggl]
\end{equation}
Therefore, using \eqref{+} and \eqref{++} we get that 
\begin{equation}
\nabla X_a \nabla \o= \omega_a \frac{a \nabla |x|}{|x| log |x|} \biggl( \textrm{log} |x| \nabla \alpha + \alpha \frac{\nabla |x|}{|x|}\biggl)=\o_a \biggl[\frac{a \nabla |x| \nabla \alpha }{|x|} + \frac{a \alpha }{ |x|^2 \textrm{log} |x|}\biggl].
\end{equation}
Then we conclude that
\begin{equation}\label{nice3}
2 \nabla X_a \nabla \o= \o_a \biggl[ \frac{2 a \nabla |x| \nabla \alpha}{|x|}+ \frac{2 a \alpha}{|x|^2 log |x|} \biggl].
\end{equation}
The sum of \eqref{nice1}, \eqref{nice2} and \eqref{nice3} we get that
\begin{equation}
\Delta \o_a=\o_a \biggl[\alpha \frac{N-K-2}{|x|^2}+ \frac{\alpha^2 }{|x|^2}+ \frac{a(a-1)}{|x|^2 (\textrm{log}|x|)^2}+ a\frac{N-k-2}{|x|^2 log|x|}+\frac{2 a \nabla |x| \nabla \alpha}{|x|}+ \frac{2 a \alpha}{|x|^2 \textrm{log} |x|}+O\bigl(\textrm{log}|x| |x|^{-3/2}\bigl)\biggl].
\end{equation}
Moreover $$\alpha (N-k-2)+\alpha^2= \biggl(\displaystyle\frac{N-k-2}{2}\biggl)^2 \biggl(-q\circ\s(x)+|x|\biggl).$$
Then using \eqref{18} we can conclude that
\begin{equation}\label{fin1}
\Delta \o_a=-\biggl(\frac{N-k-2}{2}\biggl)^2 q |x|^{-2} \o_a+a(a-1) |x|^{-2}\bigl(\textrm{log} |x|\bigl)^{-2} \o_a
\end{equation}
$$+\bigl(N-k-a\bigl) |x|^{-2} \bigl(\textrm{log} |x|\bigl)^{-1} \o_a+ O\bigl(\textrm{log} |x| |x|^{-3/2}\bigl) \o_a.$$
Using the Laplace-Beltrami operator
\begin{equation}
 \Delta_g= -g^{ij} \biggl(\frac{\partial^2}{\partial x_i \partial x_j}- \Gamma_{ij}^k \frac{\partial}{\partial x_k}\biggl),
\end{equation}
and the aproximations $$\Gamma_{ij}^k(x)=O_k(|y|)$$ and $$g_{ij}(x)=\delta_{ij}+O(|y|^2),$$ it follows that
\begin{equation}
\Delta_g v_a= \Delta_{\R^N} \o_a(F(x))+ O_{ij}\bigl(\rho^2\bigl) \partial_{ij} \o_a+ O_k(\rho) \partial_k \o_a.
\end{equation}
Now using the above identity, we conclud that for $\rho_{\Sigma_k}$ small enough
\begin{equation}\label{FIN}
\Delta_g v_a=-\biggl(\frac{N-k-2}{2}\biggl)^2 q \rho^{-2} v_a+a(a-1) \rho^{-2}\bigl(\textrm{log} \rho\bigl)^{-2} v_a
\end{equation}
$$+\bigl(N-k-a\bigl) \rho^{-2} \bigl(\textrm{log} \rho\bigl)^{-1} v_a+ O\bigl((\textrm{log} \rho) \rho^{-3/2}\bigl)v_a \qquad\textrm{ in }\Sigma_r.$$
\QED
\subsection{Construction of a sub and supersolutions}
For $\l \in \R $, $\eta \in Lip(\M)$ with $\eta=0$ on $\Sigma$ and $q \in \mathcal{C}^2(\M)$, $q> 0$ in $\M$ with $\displaystyle\max_{\Sigma} q(\sigma)=1$, we define the operator
\begin{equation}\label{lll}
L_\l:= -\Delta- \biggl(\frac{N-k-2}{2}\biggl)^2 q \rho^{-2} +\l \eta \rho^{-2}.
\end{equation}
Using Lemma \ref{lem3}} and \eqref{lll} it's easy to verify that
\begin{equation}\label{Lv}
Lv_a= -a(a-1) \rho^{-2}\bigl(\textrm{log} \rho\bigl)^{-2} v_a - \bigl(N-k-a\bigl) \rho^{-2} (\textrm{log} \rho)^{-1} v_a + \l \eta \rho^{-2}v_a+ O\bigl(\textrm{log} \rho) \rho^{-3/2}\bigl)v_a.
\end{equation}
In this subsection we wish to construct a subsolution and a supersolution for the operator $L_\l$ defined above. For that we obtain the following lemmas
\begin{Lemma}\label{lemma1}
There exists $r_0 $ such that for all $r \in (0,r_0)$ and for all $\epsilon \in [0,1)$, the function
\begin{equation}\label{V}
V_\epsilon= v_{-1,q}+v_{0,q-\epsilon}
\end{equation}
satisfies
\begin{equation}\label{inequality1}
L_\l V_\epsilon \leq 0 \textrm{ in }\Sigma_r,\;\;\ \textrm{for all} \;\;\ \epsilon \in [0,1).
\end{equation}
Moreover  $V_\epsilon \in H^1(\Sigma_r)$ for any $\epsilon \in (0,1)$ and in addition
\begin{equation}\label{inequality2}
\int_{\Sigma_r} V_0^2 \rho^{-2} dv_g \geq C \int_{\Sigma} \frac{1}{\sqrt{1-q(\s)}} d\s
\end{equation}
\end{Lemma}
\proof
Using polar coordinates it's easy to see that $v_a \in H^1(\Sigma_r)$ for all $a$ such that $a<-\displaystyle\frac{1}{2}$ and that $v_{0,q-\e} \in H^1(\Sigma_r)$ for all $\e>0$. We therefore skip the proof. Now we have that for $a=-1$ and for $r$ small enough
\begin{equation}
L_\l v_{-1,q} = -2 \rho^{-2} \bigl(\textrm{log} \rho\bigl)^{-2} v_{-1,q}+\l \eta \rho^{-2} v_{-1, q} - \bigl(N-k+1\bigl) \rho^{-2} (\textrm{log} \rho)^{-1} v_{-1}+O\bigl(\rho^{-3/2}(\textrm{log} \rho) v_{-1,q}\bigl).
\end{equation}
Therefore
\begin{equation}
L_\l v_{-1,q} \leq \biggl[-2 \rho^{-2} \bigl(\textrm{log} \rho\bigl)^{-2} + C |\textrm{log}\rho| \rho^{-3/2}+ |\l| \eta \rho^{-2}\biggl] v_{-1,q} \qquad\textrm{ in } \Sigma_r.
\end{equation}
Using the fact that $\eta=0 $ on $\Sigma$ and $\eta \in Lip(M)$ we have $|\eta|< C \rho $ around $\Sigma.$
Therefore
\begin{equation}\label{ineq1}
L_\l v_{-1, q} \leq -\rho^{-2} (\textrm{log} \rho)^{-2} v_{-1,q}= \rho^{-2} (\textrm{log} \rho)^{-3} v_{0,q}\qquad \textrm{ in }\Sigma_r.
\end{equation}
Using the same arguments as above, we get that
\begin{equation}\label{ineq2}
L_\l v_{0,q-\epsilon} \leq C|\textrm{log} \rho| \rho^{-3/2} v_{0,q-\epsilon} \qquad \textrm{ in } \Sigma_r \;\;\ \forall \epsilon \in [0,1).
\end{equation}
Therefore using \eqref{ineq1} and \eqref{ineq2} we get \eqref{inequality1}. 
\eqref{V} implies that
$$
\begin{array}{ll}
\displaystyle\int_{\Sigma_r} \displaystyle\frac{V_0^2}{\rho^2} dv_g &\geq \displaystyle\int_{\Sigma_r} \frac{v^2_{0,q}}{\rho^2} dv_g\\\
&=\displaystyle\int_{\Sigma_r} \rho^{2 \alpha-2} dv_g.
\end{array}
$$
Using \eqref{A5}, we get that
$$ \displaystyle\int_{\Sigma_r} \frac{V_0^2}{\rho^2} dv_g\geq\displaystyle\int_{\bigcup_{i=1}^{N_0} \overline{\O}_i} \rho^{2\alpha-2} dv_g,$$
$$
\begin{array}{ll}
\displaystyle \int_{\Sigma_r} \frac{V_0^2}{\rho^2} dv_g &\geq \sum_{i=1}^{N_0}\displaystyle\int_{\O_i} \rho^{2\alpha-2} dv_g\\\
&=\sum_{i=1}^{N_0}\displaystyle\int_{F^{p_i}_\M(B^{N-k}(0,r)\times D_i)} \rho^{2\alpha-2} dv_g.
\end{array}
$$
Using change of variable formula we get
$$\displaystyle\int_{\Sigma_r} \frac{V_0^2}{\rho^2} dv_g \geq \sum_{i=1}^{N_0}\displaystyle\int_{B^{N-k}(0,r)\times D_i} |z|^{2\alpha(F^{p_i}_\M(x))-2} |\textrm{Jac}(F^{p_i}_\M)|(x) dx.$$
Notice that $|\textrm{Jac}(F^{p_i}_\M)|(x)$ is bounded, the function  $|z|^{-\sqrt{|z|}} $ is also bounded in a neighborhood of the ball centered at $0$. Moreover
\begin{equation}
 \alpha(F^{p_i}_\M(x))= (2+k-N)(1-\sqrt{1-q(f^{p_i}(y))+|z|})
 \end{equation}
 so that
$$
\begin{array}{ll}
\displaystyle\int_{\Sigma_r} \displaystyle \frac{V_0^2}{\rho^2} dv_g & \geq C \sum_{i=1}^{N_0} \displaystyle\int_{B^{N-k}(0,r)\times D_i} |z|^{k-N} |z|^{-(2+k-N) \sqrt{1-q(f^{p_i}(y))}} |z|^{-\sqrt{|z|}} dx\\\\
&\geq C \sum_{i=1}^{N_0} \displaystyle \int_{B^{N-k}(0,r)\times D_i} |z|^{k-N} |z|^{-(2+k-N) \sqrt{1-q(f^{p_i}(y))}} dx.
\end{array} $$
Using polar coordinates, we get 
$$ 
\begin{array}{ll}
\displaystyle\int_{\Sigma_r} \displaystyle\frac{V_0^2}{\rho^2} dv_g &\geq C \sum_{i=1}^{N_0} \displaystyle\int_{D_i} \displaystyle\int_{S^{N-k-1}} d\theta \int_{0}^r t^{N-k-1} t^{k-N} t^{(N-k-2)\sqrt{1-q(f^{p_i}(y))}} dt dy\\\\
&\geq C \sum_{i=1}^{N_0} \displaystyle\int_{D_i} \displaystyle\int_{0}^{r_{i_1}} t^{-1} t^{(N-k-2)\sqrt{1-q(f^{p_i}(y))}} |\textrm{Jac}(f^{p_i})|(y) dy.
\end{array}$$
Therefore, using the fact that $|\textrm{Jac}(f^{p_i})|(y)=1+O(r)$ and so bounded, we get the result
\begin{equation}
\begin{array}{ll}
\displaystyle\int_{\Sigma_r} \frac{V_0^2}{\rho^2} dv_g &\geq C \displaystyle\int_{\Sigma} \int_0^r t^{-1+(N-k-2) \sqrt{1-q(\s)}} dr d\s \\\\
&\geq C \displaystyle\int_{\Sigma} \frac{r^{(N-k-2) \sqrt{1-q(\s)}}}{(N-k-2) \sqrt{1-q(\s)}} d\s\\\\
&\geq C \displaystyle\int_{\Sigma}\frac{d\s}{\sqrt{1-q(\s)}}.
\end{array}
\end{equation}
This ends the proof. \QED
\begin{Lemma}\label{lemma2}
there exists $r_0$ such that for all $r \in (0,r_0)$ the function
\begin{equation}\label{Equation2}
U=v_0-v_{-1} >0 \qquad \textrm{ in }\Sigma_r
\end{equation}
and satisfies $L_\l U \geq 0$ in $\Sigma_r.$ Moreover $U\in H^1(\Sigma_r)$ provided
\end{Lemma}
\begin{equation}\label{alala}
\int_{\Sigma} \frac{d\s}{\sqrt{1-q(\s)}} < +\infty.
\end{equation}
\proof
Using \eqref{lll}, we have that
\begin{equation}
L_\l v_0=-\bigl(N-k\bigl) \rho^{-2} \bigl(\textrm{log} \rho\bigl)^{-2} v_0+ \l \eta \rho^{-2} v_0+ O\bigl((\textrm{log} \rho) \rho^{-3/2}\bigl)v_0
\end{equation}
and
\begin{equation}
-L_{v_{-1}}= 2 \rho^{-2}\bigl(\textrm{log} \rho\bigl)^{-2} v_{-1}+ (N-k+1) \rho^{-2}\bigl(\textrm{log} \rho\bigl)^{-1}  v_{-1}-\l \eta \rho^{-2} v_{-1} +O\bigl((\textrm{log} \rho) \rho^{-3/2}\bigl) v_{-1}.
\end{equation}
so that
\begin{equation}\label{1111}
L_\l v_0 \geq - |\l| \eta \rho^{-2} v_0- C|\textrm{log} \rho)| \rho^{-3/2}v_0,
\end{equation}
\begin{equation}\label{2222}
L_\l v_{-1} \geq \bigl(2 \rho^{-2} (\textrm{log} \rho)^{-2}-C |\textrm{log} \rho| \rho^{-3/2}-|\l| \eta \rho^{-2}\bigl)v_{-1}.
\end{equation}
The dominant term in the right hand sides of the two above inequalities is $2 \rho^{-2} (log \rho)^{-2}$. Therefore there exists $r_0$ small such that for all $r \in (0,r_0)$ the inequality
\begin{equation}
L_\l U \geq 0 \textrm{ in }\Sigma_r
\end{equation}
holds.
Now we prove that $U \in H^1(\Sigma_r)$ provided inequality \eqref{alala} holds.
We have that 
$$\nabla_g v_0= \nabla (\rho^\alpha)= v_0 \nabla (\alpha \textrm{log} \rho)=v_0 \biggl(\textrm{log} \rho \nabla \alpha +\alpha \frac{\nabla \rho}{\rho}\biggl).$$ Hence  $$|\nabla v_0|^2= v_0^2 \biggl[|\textrm{log} \rho \nabla \alpha|^2+\alpha^2 \displaystyle\frac{|\nabla \rho|^2}{\rho^2}+2 \alpha (\textrm{log} \rho) \displaystyle\frac{\nabla \alpha \nabla \rho }{\rho}\biggl].$$
Using the fact that $\alpha$ is of class $\mathcal{C}^1$ and the estimation $$2 \alpha (log \rho) \displaystyle\frac{\nabla \alpha \nabla \rho}{\rho}= O(\rho^{-2} log \rho)$$ we deduce that there exists a positive constant $C$ such that
$$|\nabla v_0|^2 \leq C v_0^2 \rho^{-2}= C \rho^{2 \alpha -2}.$$
Therefore
$$\int_{{\Sigma_r}} |\nabla v_0|^2 dv_g \leq C\int_{{\Sigma_r}} \rho^{2 \alpha-2} dv_g.$$
As in the above lemma and using polar coordinates, we get
$$
\int_{\Sigma_r} |\nabla v_0|^2 dv_g \leq C \sum_{i=1}^{N_0} \int_{D_i} \int_{S^{N-k-1}} d\theta \int_0^r t^{-1} t^{(N-k-2)\sqrt{1-q(f^{p_i}(y))}} dt dy.
$$
Also as in the above lemma 
$$\sum_{i=1}^{N_0}\displaystyle\int_{D_i} \displaystyle\frac{1}{\sqrt{1-q(f^{p_i}(y))}} dy \leq C \int_{\displaystyle\Sigma} \frac{1}{\sqrt{1-q(\s)}} d\s$$ so that
\begin{equation}
\displaystyle\int_{\Sigma_r} |\nabla v_0|^2 dv_g \leq C \int_{\displaystyle\Sigma} \frac{1}{\sqrt{1-q(\s)}} d\s.
\end{equation}
This ends the proof of the Lemma. \QED
\section{Proof of Theorem \ref{Theorem 1}}\label{section 3}
In this section we give a complete proof of Theorem \ref{Theorem 1}. In the first subsection we prove the existence of $\l^*$ verifying \eqref{EL1} and \eqref{EL2}. In the second and last one of this section we give the proof of the existence and non-existence result for $\l\neq \l^*$.
\subsection{Existence of $\l^*$}
\begin{Lemma}\label{lemma31}
Let $(\M,g)$ be a smooth compact Riemannian manifold of dimension $N\geq 3$ and let $\Sigma$ be a closed submanifold of dimension $1 \leq k\leq  N-2$. We assume that the weight functions $b, q$ and $\eta$ satisfy \eqref{A22} and \eqref{A33}. Then there exists $r_0 >0 $ and $C>0$ depending only on $\M, \Sigma, q, \eta$  and $b$ such that for all $r \in (0,r_0)$ the inequality
\begin{equation}\label{B1}
\int_{\Sigma_r} b |\nabla u|^2 dv_g \geq \biggl(\frac{N-k-2}{2}\biggl)^2 \int_{\Sigma_r} q \frac{|u|^2}{\rho^2} dv_g + C \int_{\Sigma_r}\frac{|u|^2}{\rho^2 (log \rho)^2}dv_g
\end{equation}
holds for all $u \in H^1(\Sigma_r)$.
\end{Lemma}
\proof
We have that $\displaystyle \frac{b}{q} \in \mathcal{C}^{2}(\M)$, there exists $C>0$ such that:
\begin{equation}
\biggl|\frac{b(p)}{q(p)}- \frac{b(\s(p))}{q(\s(p))}\biggl| < C \rho, \;\;\ \forall \;\ p \in \Sigma_r
\end{equation}
for $r$ small enough. Hence by \eqref{A33}, there exists $C'>0$ such that
\begin{equation}\label{44}
b(p)\geq q(p)-C' \rho, \;\;\ \forall p \in \Sigma_r.
\end{equation}
Let $V=v_{1/2,q}$ in $\Sigma_r$. We have that
$$\textrm{div}(b \nabla V)=b \Delta V+ \nabla p \nabla V.$$
and by lemma \eqref{lem3} we get
\begin{equation}
-\frac{\textrm{div}(b \nabla V)}{V}\geq b \displaystyle \biggl(\frac{N-k-2}{2}\biggl)^2+\frac{1}{4} b \rho^{-2} (\textrm{log} \rho)^{-2}+ O(\rho^{-3/2}|\textrm{log} \rho|) \textrm{ in } \Sigma_r.
\end{equation}
Using \eqref{44} with the above inequality we get
\begin{equation}\label{46}
-\frac{\textrm{div}(b \nabla V)}{V}\geq  \displaystyle \biggl(\frac{N-k-2}{2}\biggl)^2 q \rho^{-2}+c \rho^{-2}(\textrm{log} \rho)^{-2} \textrm{ in } \Sigma_r,
\end{equation}
where c is a positive constant depending only on $\M, \Sigma, q, \eta$ and $b$.\\
Let $u\in \mathcal{C}^{\infty}_c(\M \setminus \Sigma)$ and define $$\varphi:=\displaystyle\frac{u}{V}.$$ Then we have that
\begin{equation}
b|\nabla u|^2=b\bigl(|V \nabla \varphi| + \nabla V \nabla (V\varphi^2)\bigl).
\end{equation}
Hence by integration by parts we get that
$$
\int_{\Sigma_r} |\nabla u|^2 b dv_g = \int_{\Sigma_r} |V \nabla \varphi|^2 dv_g+ \int_{\Sigma_r} \biggl(-\frac{\textrm{div}(p \nabla V)}{V}\biggl) u^2 dv_g.
$$
Using \eqref{46}, we get
\begin{equation}
\int_{\Sigma_r} b |\nabla u|^2 dv_g \geq \biggl(\frac{N-k-2}{2}\biggl)^2 \int_{\Sigma_r} q \frac{|u|^2}{\rho^2} dv_g + C \int_{\Sigma_r}\frac{|u|^2}{\rho^2 (\textrm{log} \rho)^2}dv_g
\end{equation}
for all $u \in \mathcal{C}^{\infty}_c(\Sigma_r)$. Using the fact that $\mathcal{C}^\infty_c\left(\M\setminus \Sigma\right)$ is dense in $\mathcal{C}^\infty_c\left(\M\right)$ the proof remain valid for a general $u.$ Furtheremore $\mathcal{C}^{\infty}_c(\Sigma_r)$ is dense in  $H^1(\Sigma_r)$. This ends the proof. \QED
\begin{Lemma}\label{lemma4}
Let $(\M,g)$ be a smooth compact Riemannian manifold of dimension $N \geq 3$ and let $\Sigma$ be a closed submanifold of dimension $1\leq k \leq N$. Assume that \eqref{A22} and \eqref{A33} hold. Then there exists $\l^*=\l^*(\M, \Sigma,b,q,\eta) \in \R$ such that
$$\mu_\l(\M,\Sigma)= \biggl(\frac{N-k-2}{2}\biggl)^2, \;\;\ \forall \l \leq \l^*,$$
$$\mu_\l(\M,\Sigma)< \biggl(\frac{N-k-2}{2}\biggl)^2, \;\;\ \forall \l > \l^* .$$
\end{Lemma}
\proof
For $b=q=1$ and $\eta=\rho^2$, we define $\displaystyle\nu_\l(\M,\Sigma):= \displaystyle\mu_\l(\M, \Sigma).$ It's known that 
\begin{equation}
\nu_\l(\R^N,\R^k)= \biggl(\displaystyle\frac{N-k-2}{2}\biggl)^2.
\end{equation}
Therefore for any $\delta >0$, we can find $u_\delta \in \mathcal{C}^{\infty}_c(\R^N)$ such that
\begin{equation}
\int_{\R^N} |\nabla u_\delta|^2 dx \leq \biggl(\biggl(\displaystyle\frac{N-k-2}{2}\biggl)^2+ \d\biggl) \int_{\R^N} |z|^{-2} u_\d^2 dx,
\end{equation}
where $x=(z,y)\in \R^{N-k} \times \R^k.$
By \eqref{A2} there exists $\s_0 \in \Sigma$ such that $b(\s_0)=q(\s_0).$
For $r>0,$ we let $\rho_r >0$ such that for all $p \in B(\s_0, \rho_r)$ we have the following
\begin{equation}\label{bii}
\begin{cases}
b(p)\leq (1+r) q(\s_0)\\

q(p) \geq (1-r) q(\s_0)\\

\eta(p) \leq r.
\end{cases}
\end{equation}
Let $\epsilon_0 >0$ small such that for all $\e \in (0,\e_0)$,  $F^{\s_0}_\M(\e \textrm{supp}(u_\d)) \subset B(\s_0,\rho_r)$
and  we let $$x=\e^{\frac{2-N}{2}}  F^{-1}(p)).$$
Therefore we define $$v(p)=\e^{\frac{2-N}{2}} u_\d(\e^{-1} F^{-1}(p)).$$
It's clear that for every $\e \in (0,\e_0),$ $v \in \mathcal{C}^{\infty}_c(\M)$. By applying the change of variable formula and \eqref{bii}, we get
$$
\begin{array}{ll}
\mu_\l(\M,\Sigma)&\leq \displaystyle\frac{\displaystyle\int_\M b |\nabla v|^2 dv_g+ \lambda \displaystyle\int_\M \rho^{-2} \eta v^2 dvg}{\displaystyle\int_\M q \rho^{-2} v^2 dv_g}\\
&\leq \displaystyle\frac{1+r}{1-r}\frac{\displaystyle \int_\M |\nabla v|^2 dv_g}{\displaystyle\int_\M \rho^{-2} v^2 dv_g}+ \displaystyle\frac{|\l| r}{(1-r) q(\s_0)}\\
&\leq \displaystyle\frac{(1+r)(1+C\e)}{(1-r)(1-C\e)} \frac{\displaystyle\int_{\R^N} |\nabla u_\d|^2 dx }{\displaystyle\int_\M |z|^{-2} u_\d^2 dx}+ \displaystyle\frac{|\l| r}{(1-r) q(\s_0)}\\
&\leq \bigl(1+O(r)\bigl)(1+O(\e)) \biggl(\bigl(\displaystyle\frac{N-k-2}{2}\bigl)^2+\d\biggl)+ O(r).
\end{array}
$$
As $\e, r , \d \rightarrow 0$ respectively, we get that 
\begin{equation}
\mu_\l(\M,\Sigma) \leq \biggl(\displaystyle\frac{N-k-2}{2}\biggl)^2,\;\ \forall \l \in \R.
\end{equation}
To finish the proof of the lemma we have just to show the existence of $\bar{\lambda}\in \R$ such that 
$$\mu_{\bar\l}(\M,\Sigma)\geq \bigl(\displaystyle\frac{N-k-2}{2}\bigl)^2.$$
Indeed we let $\varphi \in \mathcal{C}^{\infty}_c(\M)$ such that 
\begin{equation}
\varphi=
\begin{cases}
1 , \;\;\ \textrm{  in  }\Sigma_r \\
0 \;\;\;\;\ \textrm{ otherwise}.
\end{cases}
\end{equation}
For $u \in H^1(\M)$, we write $$u=u \varphi +(1-\varphi) u$$ and notice that $$u \varphi \in H^1(\Sigma_r).$$
We then have that
$$
\begin{array}{ll}
\displaystyle \int_\M u^2 \rho^{-2} q dv_g & = \displaystyle\int_\M |u \varphi +(1-\varphi) u|^2 \rho^{-2} q dv_g\\\\
&=\displaystyle\int_\M |u\varphi|^2 \rho^{-2} q dv_g+ \displaystyle\int_\M |(1-\varphi)u|^2 \rho^{-2} q dv_g+2 \displaystyle\int_\M |u^2\varphi(1-\varphi)| q \rho^{-2} dv_g\\\\
&\leq \displaystyle \int_{\Sigma_r} |u \varphi|^2 \rho^{-2} q dv_g+ 3 \displaystyle\int_{\Sigma_r^c} |(1-\varphi)|u^2 \rho^{-2} q dv_g\\
\end{array}
$$
Using Lemma \ref{lemma31}, we have that
\begin{equation}
\int_{\Sigma_r} |u \varphi|^2 \rho^{-2} q dv_g \leq \biggl(\displaystyle\frac{N-k-2}{2}\biggl)^{-2} \int_\M b |\nabla u|^2 dv_g.
\end{equation}
Therefore, there exists $C>0$ such that 
$$
\displaystyle \int_\M q u^2 \rho^{-2} dv_g \leq \biggl(\displaystyle\frac{N-k-2}{2}\biggl)^{-2} \int_\M b |\nabla u|^2 dv_g+ C \int_\M \rho^{-2} u^2 \eta  dv_g, \;\;\ \forall u \in \mathcal{C}^{\infty}_c(\M).
$$
Taking $\bar{\l} =-C$ we get the result.
Since the function $\l \rightarrow \mu_\l(\M,\Sigma)$ is decreasing, we can define $\l^*$ as 
\begin{equation}
\lambda^*:=\sup \biggl\lbrace  \l \in \R : \mu_\l(\M,\Sigma)= \displaystyle\bigl(\frac{N-k-2}{2}\bigl)^2\biggl\rbrace.
\end{equation}
this ends the proof or the Lemma. \QED
\subsection{Existence and non-existence result in the case $\l \neq \l^*$}
\begin{Theorem}\label{th2}
Let $(\M,g)$ be a smooth compact Riemannian manifold of dimension $N$ and let $\Sigma$ be a closed submanifold of dimension $1\leq k \leq N-2$. We assume that the weight functions $b, q$ and $\eta$ verify \eqref{A22} and \eqref{A33}. Then $\mu_{\l}(\M,\Sigma)$ is not achieved for every $\l < \l^*.$
\end{Theorem}
\proof
We suppose by contradiction that for some $\l_1 < \l^*$ the infinimum $\mu_\l(\M,\Sigma)$ is attained at an element $u_1 \in H^1(\M \setminus \Sigma).$ We suppose that $u_1$ is normalised so that
$$\int_\M \rho^{-2} |u_1|^2 q dv_g=1$$
and 
$$\int_\M b |\nabla u_1|^2 dv_g- \lambda_1 \int_\M \rho^{-2} |u_1|^2 \eta dv_g= \biggl(\displaystyle\frac{N-k-2}{2}\biggl)^2.$$
Then for $\l_1 <\l< \l^*,$ we have that
\begin{equation}
\biggl(\displaystyle\frac{N-k-2}{2}\biggl)^2= \mu_\l(\M,\Sigma) \leq \int_\M b |\nabla u_1|^2 dv_g- \lambda \int_\M \rho^{-2} |u_1|^2 \eta dv_g < \biggl(\displaystyle\frac{N-k-2}{2}\biggl)^2,
\end{equation}
which is impossible. So for any $\l< \l^*$, $\mu_\l(\M,\Sigma)$ is not achieved. This ends the proof of the theorem. \QED
\begin{Theorem}\label{th7}
Let $(\M,g)$ be a smooth compact Riemannian manifold of dimension $N$ and let $\Sigma_k$ be a closed submanifold of dimension $1\leq k \leq N-2$. We assume that the weight functions $b, q$ and $\eta$ verify \eqref{A22} and \eqref{A33}. Then $\mu_{\l}(\M,\Sigma)$ is achieved for every $\l > \l^*.$
\end{Theorem}
\proof
A similar proof was done by Thiam in \cite{THIAM}. So we expose here a similar one.
Let $\lbrace u_n\rbrace$ be a minimizing sequence of $\mu_\l(\M,\Sigma)$ normalized so that
$$\int_{\M} \rho^{-2} u_n^2 q dv_g=1.$$
So we have that
\begin{equation}\label{3.18}
\mu_\l(\M,\Sigma)+o(1)= \int_\M b |\nabla u_n|^2 dv_g- \l \int_\M \rho^{-2} u_n^2 \eta dv_g.
\end{equation}
Thus $\lbrace u_n\rbrace$ is bounded in $H^1(\M)$. After passing to a subsequence, we may assume that there exists $u \in H^1(\M)$ such that
\begin{equation}\label{3.21}
v_n=u_n-u \rightharpoonup 0 \textrm{ in } H^1(\M),\;\ v_n \longrightarrow 0 \textrm{ in } L^2(\M), \;\ v_n \rightharpoonup 0  \textrm{ in } H^1(\M),
\end{equation}
$$\displaystyle \frac{v_n}{\rho} \sqrt{\eta} \rightarrow 0 \textrm{ in } L^2(\M) \textrm{ and } \displaystyle\frac{v_n \sqrt q}{\rho} \rightarrow 0 \textrm{ in }  L^2(\M).$$
Using \eqref{3.18} and \eqref{3.21} we obtain that
\begin{equation}\label{3.22}
\begin{array}{ll}
\mu_\l(\M,\Sigma)&= \displaystyle\int_\M b |\nabla u_n|^2 dv_g- \l \displaystyle\int_\M \rho^{-2} u_n^2 \eta dv_g+o(1)\\\\
&=\displaystyle\int_\M b|\nabla u|^2 dv_g+ \displaystyle\int_\M b |\nabla v_n|^2 dv_g- \l \displaystyle\int_\M \rho^{-2} |u|^2 \eta dv_g + o(1)
\end{array}
\end{equation}
and
\begin{equation}\label{3.23}
1=\displaystyle\int_\M \rho^{-2} u_n^2 q dv_g+o(1)= \int_\M \rho^{-2} u^2 q dv_g+ \int_\M \rho^{-2} v_n^2 q dv_g+o(1).
\end{equation}
Let $\l < \l^*$ so that
$$\int_\M b |\nabla v_n|^2 dv_g-\l\int_\M \rho^{-2} v_n^2 \eta dv_g \geq \biggl(\frac{N-k-2}{2}\biggl)^2 \int_\M \rho^{-2} v_n^2 q dv_g +o(1).$$
Hence by \eqref{3.23} and \eqref{3.21}
\begin{equation}\label{3.24}
\int_\M b |\nabla v_n|^2 dv_g \geq \displaystyle \biggl(\frac{N-k-2}{2}\biggl)^2 \biggl(1- \displaystyle \int_\M u^2 \rho^{-2} q dv_g\biggl)+o(1).
\end{equation}
By \eqref{3.22} and \eqref{3.24} we obtain that 
\begin{equation}
\displaystyle \int_\M b |\nabla u|^2 dv_g + \displaystyle \biggl(\frac{N-k-2}{2}\biggl)^2 \biggl(1-\displaystyle \int_\M u^2 \rho^{-2} q dv_g\biggl)- \l \displaystyle \int_\M \rho^{-2} u^2 \eta dv_g \geq \mu_\l(\M,\Sigma).
\end{equation}
But
\begin{equation}
\int_\M b |\nabla u|^2 dv_g -\l\int_\M \rho^{-2} u^2 \eta  dv_g \geq \mu_\l(\M,\Sigma) \int_\M u^2 \rho^{-2} q dv_g
\end{equation}
so that 
\begin{equation}
\biggl(\mu_\l(\M,\Sigma)-\biggl(\displaystyle\frac{N-k-2}{2}\biggl)^2\biggl) \biggl(\int_\M u^2 \rho^{-2} q dv_g-1\biggl) \leq 0.
\end{equation}
Since $$\mu_\l (\M,\Sigma) < \biggl(\displaystyle\frac{N-k-2}{2}\biggl)^2,$$ we get that $$1 \leq \displaystyle \int_\M u^2 \rho^{-2} q dv_g.$$
But by Fatou's Lemma $$1 \geq \displaystyle\int_\M u^2 \rho^{-2} q dv_g.$$
Therefore
\begin{equation}
1= \displaystyle \int_\M u^2 \rho^{-2} q dv_g.
\end{equation}
We can conclude that $u$ is a minimizer for $\mu_\l(\M,\Sigma)$ and $$\int_\M b|\nabla v_n|^2 dv_g \longrightarrow 0.$$ Thus $u_n \longrightarrow u$ in $H^1(\M)$ and the proof.\QED
These two above results of this section represent a complete proof of Theorem \ref{Theorem 1}. 
\section{Proof of theorem \ref{Theorem 2}}\label{Section 4}
In this section we give a complete proof of Theorem \ref{Theorem 2}. For that we have the following results
\begin{Theorem}\label{th1}
Let $(\M,g)$ be a smooth compact Riemannian manifold of dimension $N\geq 3$, $\Sigma$ be a closed submanifold of dimension $1 \leq k\leq N-2$ and $\l \geq 0.$
Assume that the weight functions $b, q$ and $\eta$ satisfy \eqref{A22} and \eqref{A33}. We supose also that $u \in H^1(\M\setminus \Sigma)\cap \C(\M\setminus \Sigma)$ is a non-negative solution satisfying
\begin{equation}\label{Lol}
-\textrm{div}(b \nabla u)- \biggl(\displaystyle \frac{N-k-2}{2}\biggl)^2 q \rho^{-2} u \geq -\l \eta \rho^{-2} u \;\;\ \textrm{in} \;\;\ \M.
\end{equation}
Moreover if 
\begin{equation}\label{infinie}
\displaystyle \int_{\Sigma} \frac{1}{\sqrt{1-q(\s)/b(\s)}} d\s= +\infty
\end{equation}
then $u\equiv0.$
\end{Theorem}
\proof
We assume by contradiction that u does not vanish identically near $\Sigma$ and satisfies \eqref{Lol}. Therefore by standard regularity and the maximum principle, see \cite{GT}, u is smooth and positive in $\Sigma_r$ for some $r>0$ small. 
Let $\overline{u}:= \sqrt{b} u$ and then
$$
\begin{array}{ll}
\Delta \overline{u}&=\Delta (\sqrt{b} u)\\\
&= \sqrt{b} \Delta u + u\Delta \sqrt{b} + 2 \nabla u\displaystyle\frac{\nabla b}{2 \sqrt{b}}\\\\\
&=\sqrt{b} \Delta u + u \sqrt{b}\biggl[\displaystyle\frac{1}{4}|\nabla \textrm{log} b |^2+\displaystyle \frac{1}{2}\Delta \textrm{log} b\biggl]+ \displaystyle\frac{\nabla b \nabla u}{\sqrt{b}}\\\\\
&=\displaystyle\frac{1}{\sqrt{b}} \bigl( b \Delta u + \nabla b \nabla u\bigl) + u \sqrt{b} \biggl( \displaystyle\frac{|\nabla b|^2}{4 b^2}-\displaystyle\frac{|\nabla b|^2}{2 b^2}+\displaystyle\frac{\Delta b}{2b} \biggl)\\\\\
&=\displaystyle\frac{1}{\sqrt{b}} \textrm{div}(b \nabla u)+u \sqrt{b} \biggl( \displaystyle\frac{|\nabla b|^2}{4 b^2}-\displaystyle\frac{|\nabla b|^2}{2 b^2}+\displaystyle\frac{\Delta b}{2b} \biggl).\\\\\
\end{array}
$$
Therefore using \eqref{Lol}, we get that
\begin{equation}\label{Lol2}
-\Delta \overline{u}- \biggl(\displaystyle \frac{N-k-2}{2}\biggl)^2 \displaystyle\frac{q}{b} \rho^{-2} \overline{u} \geq -\l \displaystyle \frac{\eta}{b} \rho^{-2}\overline{u}+ \biggl(\frac{\Delta b}{2b}+ \frac{|\nabla b|^2}{4 b^2}\biggl)\overline{u} \;\;\ \textrm{in} \;\;\ \M.
\end{equation}
Since $b \in \mathcal{C}^2(\M)$ and $b >0$ in $\M$, the result is the same as in the case $b\equiv 1$ and $q/b$ replaced by $q.$ See Brezis-Marcus  \cite{BM} or Fall-Mahmoudhi \cite{FMa}.
So withoout lost of generality, we suppose that $b \equiv 1$ and consider the function $V_\e \in H^1(\Sigma_r)$ given by Lemma \ref{lemma1} which satisfies
\begin{equation}\label{Subsolution}
L_\l V_\e \leq 0  \textrm{ in }\Sigma_r,\textrm{ for all } \e \in (0,1).
\end{equation}
According to \eqref{Subsolution} and \eqref{Lol}, we let $R> 0$ such that 
\begin{equation}
R V_\e \leq u \textrm{ on } \partial\Sigma_r
\end{equation}
and define $$W_\e=R V_\e-u$$ so that $W^{+}_\e \in H^{1}(\Sigma_r)$.
Moreover by \eqref{Lol} and \eqref{Subsolution} we get that
\begin{equation}
L_\l W_\e \leq 0  \textrm{ in }\Sigma_r, \;\;\ \forall \e \in (0,1).
\end{equation}
Multiplying the above inequality by $W_\e^{+}$ and integrating by pats we get 
$$
\int_{\Sigma_r} |\nabla W_{\e}^{+}|^2 dv_g -\biggl(\displaystyle \frac{N-k-2}{2}\biggl)^2 \int_{\Sigma_r} \rho^{-2} q |W_\e^{+}|^2 dv_g + \l \int_{\Sigma_r} \eta \rho^{-2} |W_\e^{+}|^2 dv_g \leq 0.
$$
Then Lemma \ref{lemma31} implies that $W_\e^{+}=0$ in $\Sigma_r$ provided $r$ small enough because of the fact that $|\eta|\leq C \rho$ near $\Sigma$. Therefore $u \geq  R V_\e$ for every $\e \in (0,1)$. In particular $u \geq R V_0.$
Hence by Lemma \ref{lemma1}, we have that
\begin{equation}
\infty > \int_{\Sigma_r} u^2 \rho^{-2} dv_g \geq R^2 \int_{\Sigma_r} V_0^2 \rho^{-2} dv_g \geq \int_{\Sigma} \displaystyle \frac{d \s}{\sqrt{1-q(\s)}}.
\end{equation}
This is impossible because of \eqref{infinie}. Therefore $u\equiv 0$ in $\Sigma_r$ and by the maximum principle $u\equiv 0$ in $\M$. This ends the proof of the theorem.
\begin{Theorem}\label{th3}
Let $(\M,g)$ be a smooth compact Riemannian manifold of dimension $N$ and let $\Sigma$ be a closed submanifold of dimension $1\leq k \leq N-2$. We assume that the weight functions $b, q$ and $\eta$ verify \eqref{A22} and \eqref{A33}. If
\begin{equation}
\int_{\Sigma} \frac{1}{\sqrt{1-b(\s)/q(\s)}} d\s < \infty
\end{equation}
then $\mu_{\l^*}=\mu_{\l^*}(\M,\Sigma)$ is achieved.
\end{Theorem}
\proof
As in Theorem \ref{th1} we supose without any loss of generality that $b\equiv 1.$ Let $\lbrace \l_n\rbrace$ be a sequence of real numbers decreasing to $\l^*$. This means that $\l_n > \l^*$ for all $n \in \N$. By Theorem \ref{th2} , there exists $u_n \in H^1(\M)$ such that for all $n \in \N$
\begin{equation}\label{4.9}
-\Delta_g u_n -\mu_{\l_n}(\M) \rho^{-2} q u_n= -\l_n \rho^{-2} \eta u_n  \qquad\textrm{ in } \M.
\end{equation}
Recall that for $u_n \in H^1(\M)$, $|u_n| \in H^1(\M) $ and $|\nabla u_n|=|\nabla |u_n||.$ See for instance books  \cite{OEF} and \cite{EH} for mor details. 
Therefore we suppose that $u_n \geq 0$ in $\M$ and $||\rho^{-1} u_n||^2_2=1.$ Hence $$u_n \rightharpoonup u \qquad \textrm{ in } \;\ H^1(\M \setminus \Sigma) \;\;\ \textrm{ and }\;\;\ u_n \longrightarrow u \qquad\textrm{ in }\;\ L^2(\M).$$
We have that
\begin{equation}
\Delta_g u_n +\bigl(\mu_{\l_n} \rho^{-2} q - \l_n \rho^{-2} \eta \bigl)u_n=0 \qquad\textrm{ in } \M.
\end{equation}
We want to show that there exists $C>0$ such that
\begin{equation}\label{56}
\forall n \in \N, \;\;\ u_n \leq C U \textrm{ in } \Sigma_r.
\end{equation}
Indeed we can choose $C>0$ such that $$\forall n \in \N , \;\ v_n:=u_n-CU\leq 0 \textrm{ on } \partial\Sigma_r.$$ 
It's clear that $v_n^+ \in H^1(\Sigma_r).$ Hence 
\begin{equation}
L_{\l_n} v_n \leq -C\bigl(\mu_{\l^*}- \mu_n\bigl) q U -C \bigl(\l^* -\l_n\bigl) \eta U \leq 0 \textrm{ in }\Sigma_r.
\end{equation}
mulplying the above inequality by $v_n^+$ and integrating by parts, we get that
\begin{equation}
\displaystyle \int_{\Sigma_r} |\nabla v_n^+|^2 dv_g -\mu_{\l_n}\displaystyle \int_{\Sigma_r} \rho^{-2} q |v_n^+|^2 dv_g+ \l_n \displaystyle\int_{\Sigma_r} \eta \rho^{-2} |v_n^+|^2 dv_g \leq 0.
\end{equation}
But Lemma \ref{lemma31} gives that
\begin{equation}
C \displaystyle\int_{\Sigma_r} \rho^{-2}\bigl(log \rho\big)^{-2} |v_n^+|^2 dv_g + \l_n \displaystyle \int_{\Sigma_r} \eta \rho^{-2} |v_n^+|^2 dv_g \leq 0.
\end{equation}
Moreover $|\eta| < C\rho$ in $\Sigma_r$ and $\l_n\searrow \l^*$ so bounded. Therefore there exists $r_0 >0 $ indenpend of $n$ such that
$v_n^+\equiv 0 \textrm{ in } \Sigma_{r_0}$. Thus we obtain \eqref{56}.
By the dominated convergence theorem, the fact that $u_n \longrightarrow u$ in $L^2(\M)$ and \eqref{56} that $$\rho^{-1} u_n \longrightarrow \rho^{-1} u \qquad\textrm{ in } L^2(\M).$$
But \begin{equation}
1= \displaystyle \int_\M |\nabla u_n|^2 dv_g +o(1)= \mu_{\l_n} \displaystyle \int_\M \rho^{-2} q u_n^2 + \l_n \displaystyle \int_\M \rho^{-2} \eta u_n^2 dv_g+o(1),
\end{equation}
taking the limit, we have
\begin{equation}
1= \mu_{\l^*} \displaystyle \int_\M \rho^{-2} q u^2 + \l^* \displaystyle \int_\M \rho^{-2} \eta u^2 dv_g.
\end{equation}
Hence $u \neq 0$ and it's a minimizer for $\mu_{\l^*}.$
\subsection*{Proof of theorem \ref{Theorem 2}}
For the proof of this theorem the "if" part is given by Theorem \ref{th3} and the "only if" part is done in Theorem \ref{th1}. \QED

\end{document}